\documentclass{amsart}

\usepackage{graphicx}
\usepackage{amsmath}
\usepackage{amsthm}
\usepackage{amsfonts}
\usepackage{amssymb, amscd}

\usepackage{color}

\usepackage[all]{xy}

\newtheorem{thm}{Theorem}[section]

\newtheorem{lem}[thm]{Lemma}
\newtheorem{proposition}[thm]{Proposition}
\newtheorem{example}[thm]{Example}
\theoremstyle{definition}
\newtheorem{definition}[thm]{Definition}
\theoremstyle{remark}
\newtheorem{remark}[thm]{Remark}
\numberwithin{equation}{section}

\newcommand{\K}{\mathbb K}

\newcommand{\dl}{\displaystyle}

\begin{document}

\title[  Representations and Cohomology of n-ary multiplicative  Hom-Nambu-Lie algebras]
{Representations and Cohomology of n-ary multiplicative  Hom-Nambu-Lie algebras}%
\author{F. AMMAR, S. MABROUK and  A. MAKHLOUF}%
\address{Abdenacer Makhlouf, Universit\'{e} de Haute Alsace,  Laboratoire de Math\'{e}matiques, Informatique et Applications,
4, rue des Fr\`{e}res Lumi\`{e}re F-68093 Mulhouse, France}%
\email{Abdenacer.Makhlouf@uha.fr}
\address{Faouzi Ammar and  Sami  Mabrouk, Universit\'{e} de Sfax,  Facult\'{e} des Sciences, Sfax Tunisia}%
\email{Faouzi.Ammar@fss.rnu.tn, \ Mabrouksami00@yahoo.fr}

\thanks {
}

%
\begin{abstract}
The aim of this paper is to provide cohomologies  of $n$-ary Hom-Nambu-Lie algebras governing  central extensions and one parameter formal deformations.  We  generalize to $n$-ary algebras the notions of derivations and representation introduced by Sheng for Hom-Lie algebras. Also we show that a cohomology of $n$-ary Hom-Nambu-Lie  algebras could be derived from the cohomology of  Hom-Leibniz  algebras.
\end{abstract}
\maketitle

\section*{Introduction}
The first instances of $n$-ary algebras in Physics appeared with a generalization of the Hamiltonian mechanics proposed in 1973 by Nambu \cite{Nam}. More recent motivation comes from  string theory and M-branes involving
naturally an algebra with ternary operation called Bagger-Lambert algebra  which give  impulse to a significant development. It was used in \cite{bag} as one of the main ingredients in the construction of a new type of supersymmetric gauge theory that is consistent with all the symmetries
expected of a multiple M2-brane theory: 16 supersymmetries, conformal invariance, and an SO(8) R-symmetry that acts on the eight transverse scalars. On the other hand in the study of supergravity solutions describing M2-branes ending on M5-branes, the Lie algebra appearing in the original Nahm equations has to be replaced with a generalization involving ternary bracket in the lifted Nahm equations, see \cite{basu}. For other applications in Physics see \cite{ker1}, \cite{ker2}, \cite{ker3}.

The algebraic formulation of  Nambu mechanics is du to Takhtajan \cite{Takhtajan0,Takhtajan1} while  the abstract definition of $n$-ary Nambu
algebras or $n$-ary Nambu-Lie algebras (when the bracket is skew
symmetric) was given by Filippov in 1985 see \cite{Fil}. The Leibniz $n$-ary algebras
were introduced and studied in \cite{CassasLodayPirashvili}. For
deformation theory and cohomologies of $n$-ary  algebras of Lie
type, we refer to
\cite{makhloufAtaguema1,makhloufAtaguema2,Gautheron1,Izquirdo,Takhtajan1}.

The general Hom-algebra structures arose first in connection to
quasi-deformation and discretizations of Lie algebras of vector
fields. These quasi-deformations lead to quasi-Lie algebras, a
generalized Lie algebra structure in which the skew-symmetry and
Jacobi conditions are twisted. For   Hom-Lie algebras, Hom-associative algebras, Hom-Lie superalgebras, Hom-bialgebras ... see  \cite{AmmarMakhloufJA2010,LS1,HomNonAss,MS,HomHopf,HomAlgHomCoalg}. Generalizations of $n$-ary algebras of Lie type and associative type by twisting the identities using linear maps have been introduced in \cite{makh}.
These generalizations include $n$-ary
Hom-algebra structures generalizing the $n$-ary algebras of Lie type
such as $n$-ary Nambu algebras, $n$-ary Nambu-Lie algebras and
$n$-ary Lie algebras, and $n$-ary algebras of associative type
such as $n$-ary totally associative and $n$-ary partially
associative algebras. See also \cite{yau1,yau2,yau3}.

In the first Section of this paper we summarize  the definitions of $n$-ary Hom-Nambu (resp. Hom-Nambu-Lie) algebras and the multiplicative $n$-ary Hom-Nambu (resp. Hom-Nambu-Lie) algebras.
In Section $2$, we extend to $n$-ary algebras   the notions of   derivations and representation introduced for Hom-Lie algebras in \cite{Sheng}.
In Section $3$, we show that for an $n$-ary Hom-Nambu-Lie algebra $\mathcal{N}$, the space $\wedge^{n-1}\mathcal{N}$ carries a structure of  Hom-Leibniz algebra.
Section $4$ is dedicated to central extensions.  We provide a cohomology adapted to central extensions of $n$-ary multiplicative Hom-Nambu-Lie algebras.
In Section $5$, we provide a cohomology which is suitable for the study of  one parameter formal  deformations  of $n$-ary Hom-Nambu-Lie algebras.
In the last Section we show that the cohomology of $n$-ary Hom-Nambu-Lie algebras may be derived from the  the cohomology of Hom-Leibniz  algebras. To this end we generalize to twisted situation the process used  by Daletskii and Takhtajan \cite{Takhtajan0} for the classical case.
\section{The n-ary Hom-Nambu algebras}
  Throughout this paper, we will for simplicity of exposition assume that $\mathbb{K}$ is an algebraically closed
field of characteristic zero, even though for most of the general definitions and results in the paper this
assumption is not essential.
\subsection{Definitions}
In this section, we recall the definition of $n$-ary Hom-Nambu algebras and $n$-ary Hom-Nambu-Lie algebras, introduced in \cite{makh} by Ataguema, Makhlouf and Silvestrov. They correspond to a generalized version by twisting of $n$-ary Nambu algebras and Nambu-Lie algebras which are called Filippov algebras. We deal in this paper with
 a subclass of $n$-ary Hom-Nambu algebras  called  multiplicative $n$-ary Hom-Nambu algebras.
\begin{definition}
An \emph{$n$-ary Hom-Nambu} algebra is a triple $(\mathcal{N}, [\cdot ,..., \cdot],  \widetilde{\alpha} )$ consisting of a vector space  $\mathcal{N}$, an
$n$-linear map $[\cdot ,..., \cdot ] :  \mathcal{N}^{ n}\longrightarrow \mathcal{N}$ and a family
$\widetilde{\alpha}=(\alpha_i)_{1\leq i\leq n-1}$ of  linear maps $ \alpha_i:\ \ \mathcal{N}\longrightarrow \mathcal{N}$, satisfying \\
  \begin{eqnarray}\label{NambuIdentity}
  && \big[\alpha_1(x_1),....,\alpha_{n-1}(x_{n-1}),[y_1,....,y_{n}]\big]= \\ \nonumber
&& \sum_{i=1}^{n}\big[\alpha_1(y_1),....,\alpha_{i-1}(y_{i-1}),[x_1,....,x_{n-1},y_i]
  ,\alpha_i(y_{i+1}),...,\alpha_{n-1}(y_n)\big],
  \end{eqnarray}
  for all $(x_1,..., x_{n-1})\in \mathcal{N}^{ n-1}$, $(y_1,...,  y_n)\in \mathcal{N}^{ n}.$\\
  The identity $(1.1)$ is called \emph{Hom-Nambu identity}.
  \end{definition}
Let
$x=(x_1,\ldots,x_{n-1})\in \mathcal{N}^{n-1}$, $\widetilde{\alpha}
(x)=(\alpha_1(x_1),\ldots,\alpha_{n-1}(x_{n-1}))\in \mathcal{N}^{n-1}$ and
$y\in \mathcal{N}$. We define an adjoint map  $ad(x)$ as  a linear map on $\mathcal{N}$,
such that
\begin{equation}\label{adjointMapNaire}
ad(x)(y)=[x_{1},\cdots,x_{n-1},y].
\end{equation}

Then the Hom-Nambu identity \eqref{NambuIdentity} may be written in terms of adjoint map as
\begin{equation*}
ad(\widetilde{\alpha} (x))( [x_{n},...,x_{2n-1}])=
\sum_{i=n}^{2n-1}{[\alpha_1(x_{n}),...,\alpha_{i-n}(x_{i-1}),
ad(x)(x_{i}), \alpha_{i-n+1}(x_{i+1}) ...,\alpha_{n-1}(x_{2n-1})].}
\end{equation*}

\begin{remark}
When the maps $(\alpha_i)_{1\leq i\leq n-1}$ are all identity maps, one recovers the classical $n$-ary Nambu algebras. The Hom-Nambu Identity \eqref{NambuIdentity}, for $n=2$,  corresponds to Hom-Jacobi identity (see \cite{MS}), which reduces to Jacobi identity when $\alpha_1=id$.
\end{remark}

Let $(\mathcal{N},[\cdot,\dots,\cdot],\widetilde{\alpha})$ and
$(\mathcal{N}',[\cdot,\dots,\cdot]',\widetilde{\alpha}')$ be two $n$-ary Hom-Nambu
algebras  where $\widetilde{\alpha}=(\alpha_{i})_{i=1,\cdots,n-1}$ and
$\widetilde{\alpha}'=(\alpha'_{i})_{i=1,\cdots,n-1}$. A linear map $f:
\mathcal{N}\rightarrow \mathcal{N}'$ is an  $n$-ary Hom-Nambu algebras \emph{morphism}  if it satisfies
\begin{eqnarray*}f ([x_{1},\cdots,x_{n}])&=&
[f (x_{1}),\cdots,f (x_{n})]'\\
f \circ \alpha_i&=&\alpha'_i\circ f \quad \forall i=1,n-1.
\end{eqnarray*}

\begin{definition}
An $n$-ary Hom-Nambu algebra $(\mathcal{N}, [\cdot ,..., \cdot],  \widetilde{ \alpha} )$ where  $\widetilde{\alpha}=(\alpha_i)_{1\leq i\leq n-1}$
is called \emph{$n$-ary Hom-Nambu-Lie} algebra if the bracket is skew-symmetric that is
\begin{equation}
[x_{\sigma(1)},..,x_{\sigma(n)}]=Sgn(\sigma)[x_1,..,x_n],\ \ \forall \sigma\in \mathcal{S}_n
\ \ \textrm{and}\ \ \forall x_1,...,x_n\in \mathcal{N}.
\end{equation}
where $\mathcal{S}_n$ stands for the permutation group of $n$ elements.
\end{definition}
In the sequel we deal with a particular class of $n$-ary Hom-Nambu-Lie algebras which we call $n$-ary multiplicative Hom-Nambu-Lie algebras.
\begin{definition}
An \emph{$n$-ary multiplicative Hom-Nambu algebra }
(resp. \emph{$n$-ary multiplicative Hom-Nambu-Lie algebra}) is an $n$-ary Hom-Nambu algebra  (resp. $n$-ary Hom-Nambu-Lie algebra) $(\mathcal{N}, [\cdot ,..., \cdot],  \widetilde{ \alpha})$ with  $\widetilde{\alpha}=(\alpha_i)_{1\leq i\leq n-1}$
where  $\alpha_1=...=\alpha_{n-1}=\alpha$  and satisfying
\begin{equation}
\alpha([x_1,..,x_n])=[\alpha(x_1),..,\alpha(x_n)],\ \  \forall x_1,...,x_n\in \mathcal{N}.
\end{equation}
For simplicity, we will denote the $n$-ary multiplicative Hom-Nambu algebra as $(\mathcal{N}, [\cdot ,..., \cdot ],  \alpha)$ where $\alpha :\mathcal{N}\rightarrow \mathcal{N}$ is a linear map. Also by misuse of language an element  $x\in \mathcal{N}^n$ refers $x=(x_1,..,x_{n})$ and  $\alpha(x)$ denotes $(\alpha (x_1),...,\alpha (x_n))$.
\end{definition}

The following theorem gives a way  to construct $n$-ary multiplicative Hom-Nambu algebras (resp. Hom-Nambu-Lie  algebras) starting from  classical $n$-ary Nambu algebras (resp. Nambu-Lie  algebras) and  algebra endomorphisms.

\begin{thm}\cite{makh} Let $(\mathcal{N},[\cdot,...,\cdot])$ be an $n$-ary Nambu algebra
(resp. $n$-ary Nambu-Lie algebra) and let $\rho:\mathcal{N}\rightarrow \mathcal{N}$ be an n-ary Nambu
(resp.  Nambu-Lie) algebra endomorphism.
Then $(\mathcal{N},\rho\circ[\cdot ,...,\cdot ],\rho)$ is a $n$-ary multiplicative Hom-Nambu algebra
(resp. $n$-ary multiplicative Hom-Nambu-Lie algebra).
\end{thm}
\section{Representations of Hom-Nambu-Lie algebras}
In this section we extend the representation theory of  Hom-Lie algebras introduced in \cite{Sheng} and \cite{BenayadiMakhlouf}  to the $n$-ary case. We denote by $End(\mathcal{N})$ the linear group of operators on the $\K$-vector space $\mathcal{N}$. Sometimes it is considered as a Lie algebra with the commutator brackets.
\subsection{Derivations of $n$-ary Hom-Nambu-Lie algebras}
Let $(\mathcal{N}, [\cdot  ,..., \cdot ],  \alpha )$ be an $n$-ary multiplicative Hom-Nambu-Lie algebra. We
denote by $\alpha^k$ the $k$-times composition of $\alpha$
(i.e.  $\alpha^k=\alpha\circ...\circ\alpha$ $k$-times).
In particular $\alpha^{-1}=0$ and $\alpha^0=id$.
\begin{definition}
For any $k\geq1$, we call $D\in End(\mathcal{N})$ an $\alpha^k$-\emph{derivation }of the
$n$-ary multiplicative Hom-Nambu-Lie algebra  $(\mathcal{N}, [\cdot ,...,\cdot],  \alpha )$ if
\begin{equation}\label{alphaKderiv1}[D,\alpha]=0\ \ (\textrm{i.e.}\ \ D\circ\alpha=\alpha\circ D),\end{equation}
and
\begin{equation}\label{alphaKderiv2}
D[x_1,...,x_n]=\sum_{i=1}^n[\alpha^k(x_1),...,\alpha^k(x_{i-1}),D(x_i),\alpha^k(x_{i+1}),...,\alpha^k(x_n)],
\end{equation}
We denote by $Der_{\alpha^k}(\mathcal{N})$ the set of $\alpha^k$-derivations of
the $n$-ary multiplicative Hom-Nambu-Lie  algebra $\mathcal{N}$.
\end{definition}
For $x=(x_1,...,x_{n-1})\in \mathcal{N}^{ n-1}$ satisfying $\alpha(x)=x$ and $k\geq 1$,
we define the map $ad_k(x)\in End(\mathcal{N})$ by
\begin{equation}\label{ad_k(u)}
ad_k(x)(y)=[x_1,...,x_{n-1},\alpha^k(y)]\ \ \forall y\in \mathcal{N}.
\end{equation}
Then
\begin{lem}
The map $ad_k(x)$ is an $\alpha^{k+1}$-derivation, that we call inner $\alpha^{k+1}$-derivation.
\end{lem}
We denote by $Inn_{\alpha^k}(\mathcal{N})$ the $\K$-vector space generated by all  inner $\alpha^{k+1}$-derivations.
For any $D\in Der_{\alpha^k}(\mathcal{N})$ and $D'\in Der_{\alpha^k}(\mathcal{N})$ we define their commutator $[D,D']$ as usual:
\begin{equation}\label{DerivationsCommutator}[D,D']=D\circ D'-D'\circ D.\end{equation}
Set $Der(\mathcal{N})=\dl\bigoplus_{k\geq -1}Der_{\alpha^k}(\mathcal{N})$ and $Inn(\mathcal{N})=\dl\bigoplus_{k\geq -1}Inn_{\alpha^k}(\mathcal{N})$.
\begin{lem}\label{2.3}
For any $D\in Der_{\alpha^k}(\mathcal{N})$ and $D'\in Der_{\alpha^{k'}}(\mathcal{N})$, where $k+k'\geq-1$, we have
$$[D,D']\in Der_{\alpha^{k+k'}}(\mathcal{N}).$$
\end{lem}
\begin{proof}Let $x_i\in \mathcal{N},\ 1\leq i\leq n$, $D\in Der_{\alpha^k}(\mathcal{N})$ and $D'\in Der_{\alpha^{k'}}(\mathcal{N})$, then
\begin{eqnarray*}
  D\circ D'([x_1,...,x_n]) &=& \sum_{i=1}^n D([\alpha^{k'}(x_1),...,D'(x_i),...,\alpha^{k'}(x_n)])\\
    &=& \sum_{i=1}^n [\alpha^{k+k'}(x_1),...,D\circ D'(x_i),...,\alpha^{k+k'}(x_n)] \\
    &+& \sum_{i<j}^n [\alpha^{k+k'}(x_1),...,\alpha^k( D'(x_i)),...\alpha^{k'}(D(x_j)),...,\alpha^{k+k'}(x_n)]\\
    &+& \sum_{i>j}^n [\alpha^{k+k'}(x_1),...,\alpha^{k'}( D(x_j)),...\alpha^k(D'(x_i)),...,\alpha^{k+k'}(x_n)].
\end{eqnarray*}
The second and the third term in $[D,D']$ are symmetrical, then
\begin{eqnarray*}
  [D, D']([x_1,...,x_n]) &=& (D\circ D'-D'\circ D)([x_1,...,x_n])\\
    &=& \sum_{i=1}^n [\alpha^{k+k'}(x_1),...,(D\circ D'-D'\circ D)(x_i),...,\alpha^{k+k'}(x_n)] \\
    &= &\sum_{i=1}^n [\alpha^{k+k'}(x_1),...,[D, D'](x_i),...,\alpha^{k+k'}(x_n)],
\end{eqnarray*}
which yield that $[D,D']\in Der_{\alpha^{k+k'}}(\mathcal{N})$.
\end{proof}
Moreover we have:
\begin{proposition}
The pair $(Der(\mathcal{N}),[\cdot ,\cdot ])$, where the bracket is the usual commutator,  defines a Lie algebra and $Inn(V)$ constitutes  an  ideal  of  it.
\end{proposition}
\begin{proof}$(Der(\mathcal{N}),[\cdot ,\cdot ])$ is a Lie algebra by using Lemma \ref{2.3}.
We show that $Inn(V)$ is an ideal. Let $ad_k(x)=[x_1,...,x_{n-1},\alpha^{k-1}(\cdot)]$ be an inner $\alpha^k$-derivation on $\mathcal{N}$ and $D\in Der_{\alpha^{k'}}(\mathcal{N})$ for $k\geq -1$ and $k'\geq -1$ where $k+k'\geq -1$. Then $$[D,ad_k(x)]\in Der_{\alpha^{k+k'}}(\mathcal{N}) $$
and for any $y\in \mathcal{N}$
\begin{eqnarray*}
  [D,ad_k(x)](y) &=& D([x_1,...,x_{n-1},\alpha^{k-1}(y)])-[x_1,...,x_{n-1},\alpha^{k-1}(D(y))],\\
    &=& D([\alpha^{k}(x_1),...,\alpha^{k}(x_{n-1}),\alpha^{k-1}(y)])-[\alpha^{k+k'}(x_1),...,\alpha^{k+k'}(x_{n-1}),\alpha^{k-1}(D(y))], \\
    &=& \sum_{i\leq n-1}[\alpha^{k+k'}(x_1),...,D(\alpha^k(x_i)),...,\alpha^{k+k'}(x_{n-1}),\alpha^{k+k'-1}(y)], \\
    &=& \sum_{i\leq n-1}[x_1,...,D(x_i),...,x_{n-1},\alpha^{k+k'-1}(y)], \\
    &=&\sum_{i\leq n-1}ad_{k+k'}(x_1\wedge...\wedge D(x_i)\wedge...\wedge x_{n-1})(y).
\end{eqnarray*}
Therefore  $[D,ad_k(x)]\in Inn_{\alpha^{k+k'}}(V)$.
\end{proof}
\subsection{Representations of $n$-ary  Hom-Nambu-Lie algebras}
In this section we introduce and study the representations of $n$-ary multiplicative  Hom-Nambu-Lie algebras.
\begin{definition}
A representation of an $n$-ary multiplicative Hom-Nambu-Lie algebra $(\mathcal{N},[\cdot ,...,\cdot ],\alpha)$
on a vector space $\mathcal{N}$ is a skew-symmetric multilinear map $\rho:\mathcal{N}^{ n-1}\longrightarrow End(\mathcal{N})$,
satisfying  for $x,y\in \mathcal{N}^{n-1}$ the identity
\begin{equation}\label{RepIdentity1}
\rho(\alpha(x))\circ\rho(y)-\rho(\alpha(y))\circ\rho(x)=\sum_{i=1}^{n-1}\rho(\alpha (x_1),...,ad(y)(x_i),...,\alpha (x_{n-1}))\circ \nu
\end{equation}

where $\nu$ is a linear map.
\end{definition}
Two representations $\rho$ and $\rho'$ on $\mathcal{N}$ are \emph{equivalent} if there exists $f:\mathcal{N} \rightarrow \mathcal{N} $ an isomorphism of vector space such that $f(x\cdot y)=x\cdot ' f(y)$ where $x\cdot y=\rho(x)(y)$ and $x\cdot' y=\rho'(x)(y)$ for $x\in \mathcal{N}^{n-1}$ and $y\in \mathcal{N}.$
\begin{example}Let $(\mathcal{N}, [\cdot  ,..., \cdot ],  \alpha )$ be an $n$-ary multiplicative Hom-Nambu-Lie algebra. The map $\rm ad$  defined in \eqref{adjointMapNaire} is a representation, where the operator  $\nu$ is the twist map $\alpha$. The identity \eqref{RepIdentity1} is equivalent to Hom-Nambu identity. It is called the adjoint representation.
\end{example}
\section{From $n$-ary Hom-Nambu-Lie algebra to Hom-Leibniz algebra}
In the context of  Hom-Lie algebras one
gets the class of Hom-Leibniz algebras (see \cite{MS}). Following the standard
Loday's conventions for Leibniz algebras,
a Hom-Leibniz algebra  is a triple $(V, [\cdot, \cdot], \alpha)$
consisting of a vector  space $V$, a bilinear map $[\cdot, \cdot]:
V\times V \rightarrow V$ and a linear map $\alpha: V \rightarrow
V$ with respect to $[\cdot, \cdot]$ satisfying
\begin{equation} \label{Leibnizalgident}
[\alpha(x),[y,z]]=[[x,y],\alpha(z)]+[\alpha (y),[x,z]]
\end{equation}

Let $(\mathcal{N},[\cdot ,...,\cdot ],\alpha)$ be a $n$-ary multiplicative Hom-Nambu-Lie  algebras, we define\\

$\bullet$ a linear map  $L:\wedge^{n-1}\mathcal{N}\longrightarrow End(\mathcal{N})$ by
\begin{equation}\label{adj}L(x)\cdot z=[x_1,...,x_{n-1},z],\end{equation}
for all $x=x_1\wedge...\wedge x_{ n-1}\in\wedge^{n-1}\mathcal{N}, \ z\in \mathcal{N}$
and extending it linearly to all $\wedge^{n-1}\mathcal{N}$.  Notice that $L(x)\cdot z=ad(x)(z)$.\\

$\bullet$  a linear map $\tilde{\alpha}:\wedge^{n-1}\mathcal{N}\longrightarrow\wedge^{n-1}\mathcal{N}$ by
\begin{equation}\tilde{\alpha}(x)=\alpha(x_1)\wedge...\wedge\alpha(x_{n-1})\,\end{equation}

for all $x=x_1\wedge...\wedge x_{ n-1}\in\wedge^{n-1}\mathcal{N}$,\\

$\bullet$ a bilinear map $[\ ,\ ]_{\alpha}:\wedge^{n-1}\mathcal{N}\times\wedge^{n-1}\mathcal{N}\longrightarrow\wedge^{n-1}\mathcal{N}$  by
\begin{equation}\label{brackLei}[x ,y]_{\alpha}=L(x)\bullet_{\alpha}y=\sum_{i=0}^{n-1}\big(\alpha(y_1),...,L(x)\cdot y_i,...,\alpha(y_{n-1})\big),\end{equation}
for all $x=x_1\wedge...\wedge x_{ n-1}\in\wedge^{n-1}\mathcal{N},\ y=y_1\wedge...\wedge y_{ n-1}\in\wedge^{n-1}\mathcal{N}$\\

We denote by $\mathcal{L}(\mathcal{N})$ the space $\wedge^{n-1}\mathcal{N}$ and we call it  the fundamental set.
\begin{lem}\label{3.1}
The map $L$ satisfies
\begin{equation} L([x ,y ]_{\alpha})\cdot \alpha(z)=L(\alpha(x))\cdot \big(L(y)\cdot z\big)-L(\alpha(y))\cdot \big(L(x)\cdot z\big)\end{equation}
for all $x,\ y\in \mathcal{L}(\mathcal{N}),\ z\in \mathcal{N}$
\end{lem}
\begin{proposition}\label{HomLeibOfHomNambu}The triple $\big(\mathcal{L}(\mathcal{N}),\ [\ ,\ ]_{\alpha},\ \alpha\big)$ is a Hom-Leibniz algebra.
\end{proposition}
\begin{proof} Let $x=x_1\wedge...\wedge x_{ n-1},\ y=y_1\wedge...\wedge y_{n-1} $ and $z=z_1\wedge...\wedge z_{n-1}\in \mathcal{L}(\mathcal{N})$, the Leibniz identity \eqref{Leibnizalgident} can be written
\begin{equation}\label{brackLei2}\big[[x ,y ]_{\alpha} ,\alpha(z)\big]_{\alpha}=[\alpha(x) ,[y ,z ]_{\alpha} \big]_{\alpha}-[\alpha(y) ,[x ,z ]_{\alpha}\big]_{\alpha}\end{equation}
and  equivalently
\begin{equation}\label{brackLei3}
\Big(L\big(L(x)\bullet_{\alpha}y\big)\bullet_{\alpha}\tilde{\alpha}(z)\Big)\cdot (v)
=\Big(L(\alpha(x))\bullet_{\alpha}\big(L(y)\bullet_{\alpha}z\big)\Big)\cdot (v)-
\Big(L(\alpha(y))\bullet_{\alpha}\big(L(x)\bullet_{\alpha}z\big)\Big)\cdot (v).
\end{equation}
Let us compute first $\Big(L(\tilde{\alpha}(x))\bullet_{\alpha}\big(L(y)\bullet_{\alpha}z\big)\Big)$. This is given by
\begin{eqnarray*}
                                  \Big(L(\alpha(x))\bullet_{\alpha}\big(L(y)\bullet_{\alpha}z\big)\Big) &=& \sum_{i=0}^{n-1}L(\alpha(x))\bullet_{\alpha}\big(\alpha(z_1),...,L(y)\cdot z_i,...,\alpha(z_{n-1})\big) \\
                                    &=& \sum_{i=0}^{n-1}\sum_{i\neq j,j=0}^{n-1}\big(\alpha^2(z_1),...,\alpha(L(x)\cdot z_j),...,\alpha(L(y)\cdot z_i)...,\alpha^2(z_{n-1})\big) \\
                                    &+& \sum_{i=0}^{n-1}\big(\alpha^2(z_1),...,L(\tilde{\alpha}(x))\cdot (L(y)\cdot z_i),...,\alpha^2(z_{n-1})\big).
                                \end{eqnarray*}
                                The right hand side of \eqref{brackLei3} is skewsymmetric in $x$, $y$; hence,
                                $$
    \Big(L(\alpha(x))\bullet_{\alpha}\big(L(y)\bullet_{\alpha}z\big)\Big)-
\Big(L(\alpha(y))\bullet_{\alpha}\big(L(x)\bullet_{\alpha}z\big)\Big)=$$\begin{equation}
\sum_{i=0}^{n-1}  (\alpha^2(z_1),...,\{L(\alpha(x))\cdot (L(y)\cdot z_i)
  - L(\alpha(y))\cdot (L(x).z_i)\},...,\alpha^2(z_{n-1})\big).
\end{equation}
In the other hand, using Definition \eqref{brackLei}, we find
$$\Big(L\big(L(x)\bullet_{\alpha}y\big)\bullet_{\alpha}\tilde{\alpha}(z)\Big)
=$$$$\sum_{i=0}^{n-1}\sum_{j=0}^{n-1}\big(\alpha^2(z_1),...,\alpha^2(z_{i-1}),[\alpha(y_1),...,
L(x)\cdot y_j,...,\alpha(y_{n-1}),\alpha(z_i)],\alpha^2(z_{i+1}),...,\alpha^2(z_{n-1})\big)$$ \begin{equation}=\sum_{i=0}^{n-1}\big(\alpha^2(z_1),...,\alpha^2(z_{i-1}),[x ,y ]_{\alpha}\cdot \alpha(z_i),\alpha^2(z_{i+1}),...,\alpha^2(z_{n-1})\big).\end{equation}
Using Lemma \ref{3.1}, the proof is completed.
\end{proof}
\begin{remark}
We obtain a similar result if we consider the space $T\mathcal{N}=\otimes^n \mathcal{N}$ instead of $\mathcal{L}(\mathcal{N})$.
\end{remark}
\begin{remark}For  $n=2$  the map $L:\mathcal{L}(\mathcal{N})\longrightarrow End(\mathcal{N})$ defines  a representation of $\mathcal{L}(\mathcal{N})$ on $\mathcal{N}$.

One should set $\nu =\alpha$ and check
\begin{align}\label{I}
L(\alpha(x))\cdot \alpha(z)&=\alpha(L(x)\cdot z)\\\label{II}  L([x,y]_\alpha)\cdot \alpha(z)&=L(\alpha(x))(y)\cdot z-L(\alpha(y))(x)\cdot z
\end{align}
Indeed \eqref{I} and \eqref{II} are equivalent to
\begin{align}
[\alpha(x),\alpha (y)]&=\alpha([x,y]),\\
[[x,y],\alpha(z)]&=[[\alpha(x),y],z]-[[\alpha(y),x],z].
\end{align}
According to \cite{Sheng} and \cite{BenayadiMakhlouf} it corresponds to the adjoint representation of a Hom-Lie algebra.
\end{remark}

\section{Central Extensions and Cohomology of $n$-ary Hom-Nambu-Lie algebras}
\subsection{Central extensions of  $n$-ary multiplicative Hom-Nambu-Lie  algebras}

Let $(\mathcal{N},[\cdot ,...,\cdot ],\alpha)$ be an  $n$-ary  multiplicative Hom-Nambu-Lie  algebra.
\begin{definition}We define a central extension $\tilde{\mathcal{N}}$ of $\mathcal{N}$ by
adding a new central generator $e$ and modifying the bracket as follows: for all $\tilde{x}_i=x_i+a_i e$, $ a_i\in\K$  and $  1\leq i\leq n$ we have
\begin{equation}\label{7.1}
[\tilde{x}_1,...,\tilde{x}_n]_{\widetilde{\mathcal{N}}}=[x_1,...,x_n]+\varphi(x_1,...,x_n)e,
\end{equation}
\begin{equation}\label{7.11}\beta(\widetilde{x}_i)=\alpha(x_i)+\lambda(x_i)e,\end{equation}
\begin{equation}\label{7.12} [\tilde{x}_1,...,\tilde{x}_{n-1},e]_{\widetilde{\mathcal{N}}}=0,\end{equation}
 where $\lambda:\mathcal{N}\rightarrow\mathbb{K}$ a linear map.

\end{definition}
One may think of adding more than one central generator, but this will not be needed here for the discussion.

\begin{itemize}
\item  Clearly, $\varphi$ has to be an $n$-linear and  skew-symmetric map, $\varphi\in \wedge^{n-1}\mathcal{N}^*\wedge \mathcal{N}^*$, where $\mathcal{N}^*$ is the dual
of $\mathcal{N}$. It will be identified with a $1$-cochain.
\item  The new bracket for the $\tilde{x}_i\in\widetilde{\mathcal{N}}$ has to satisfy the Hom-Nambu identity. This leads to a condition on $\varphi$ when one of the vector involved is $e$.
\item  Since $e$ is a central then the he Hom-Nambu identity has no restriction on $\lambda$.\\
For  $\tilde{x}_i=x_i+a_i e\in \tilde{\mathcal{N}}\ $, $\tilde{y}_i=y_i+b_i e\in \tilde{\mathcal{N}},\ \ 1\leq i\leq n$, we have
  \begin{eqnarray*}&& \big[\beta(\tilde{x}_1),....,\beta(\tilde{x}_{n-1}),[\tilde{y}_1,....,\tilde{y}_{n}]_{\widetilde{\mathcal{N}}}\big]_{\widetilde{\mathcal{N}}}= \\ && \sum_{i=1}^{n-1}\big[\beta(\tilde{y}_1),....,\beta(\tilde{y}_{i-1}),[\tilde{x}_1,....,\tilde{x}_{n-1},y_i]_{\widetilde{\mathcal{N}}}.
  \beta(\tilde{y}_{i+1}),...,\beta(\tilde{y}_n)\big]_{\widetilde{\mathcal{N}}},
  \end{eqnarray*}

Using \eqref{7.1} and the Hom-Nambu identity for the original Hom-Nambu-Lie algebra, one gets
  \begin{eqnarray}\label{7.2}&& \varphi\big(\alpha(x_1),....,\alpha(x_{n-1}),[y_1,....,y_{n}]\big)- \nonumber\\ && \sum_{i=1}^{n-1}\varphi\big(\alpha(y_1),....,\alpha(y_{i-1}),[x_1,....,x_{n-1},y_i],\alpha(y_{i+1}),...,\alpha(y_n)\big)=0,
  \end{eqnarray}
  \item  The previous equation, may be written as $$\delta^2\varphi(x,y,z)=0$$
  where  $x=x_1\otimes...\otimes x_{n-1}\in \mathcal{N}^{\otimes n-1},\ y=y_1\otimes...\otimes y_{ n-1}\in \mathcal{N}^{\otimes n-1},\ z=y_n\in \mathcal{N}$.\\
  We provide below the condition that characterizes $\varphi\in \wedge^{n-1}\mathcal{N}^*\wedge \mathcal{N}^*$, $\varphi:x\wedge z\rightarrow\varphi(x,z)$ as a $1$-cocycle. It is seen now why becomes natural to call $\varphi$ a $1$-cocycle (rather than a $2$-cochain, as it is in the Hom-Lie cohomology case in \cite{HomDeform}).\\
  The number of  elements of $\mathcal{L}(\mathcal{N})$ in the argument of a cochain determines its order. As we shall see shortly, an arbitrary $p$-cochain takes $p(n-1)+1$ arguments in $\mathcal{N}$.  A $0$-cochain is an element of $\mathcal{N}^*$.
\end{itemize}

\subsection{Cohomology adapted to central extensions of multiplicative Hom-Nambu-Lie algebras}
Let us now construct the cohomology complex relevant for central extensions of multiplicative Hom-Nambu-Lie  algebras. Since $\mathcal{N}$ does not act on $\varphi(x,z)$, it will be the cohomology of multiplicative Hom-Nambu-Lie  algebras for the trivial action.
\begin{definition}
We define an arbitrary $p$-cochain as an element  $\varphi\in \wedge^{n-1}\mathcal{N}^*\otimes...\otimes\wedge^{n-1}\mathcal{N}^*\wedge \mathcal{N}^*$,
\begin{eqnarray*}
  \varphi:\mathcal{L}(\mathcal{N})\otimes...\otimes \mathcal{L}(\mathcal{N})\wedge \mathcal{N} &\longrightarrow& \mathbb{K} \\
  (x_1,..,x_p,z)\ \ \ \ &\longmapsto& \varphi(x_1,..,x_p,z)
\end{eqnarray*}

We denote the set of  $p$-cochains  with values in $\K$ by $C^p(\mathcal{N},\mathbb{K})$.
\end{definition}
Condition \eqref{7.2} guarantees the consistency of $\varphi$ according to   \eqref{7.1}  with the Hom-Nambu identity \eqref{NambuIdentity}. Then
\begin{equation}\label{7.3}\delta^2\varphi(x,y,z)=\varphi\big(\alpha(x),L(y)\cdot z\big)-\varphi\big(\alpha(y),L(x)\cdot z\big)-
  \varphi\big([x,y]_{\alpha},\alpha(z)\big)=0,\end{equation}
  where $L(x)\cdot z$ and $[x,y]_\alpha$ are defined in \eqref{adj} and \eqref{brackLei}. It is now straightforward to extend \eqref{7.3}  to a whole cohomology complex; $\delta^p\varphi$ will be a $(p+1)$-cochain taking one more argument  of $\mathcal{L}(\mathcal{N})$ than $\varphi$. This is done by means of the following
  \begin{definition}
  Let $\varphi\in C^p(\mathcal{N},\mathbb{K})$ be a $p$-cochain on a multiplicative $n$-ary Hom-Nambu-Lie algebra $\mathcal{N}$. A  coboundary operator $\delta^p$ on arbitrary $p$-cochain is given by

  \begin{eqnarray}
    \delta^p\varphi(x_1,...,x_{p+1},z) &=& \sum_{1\leq i<j}^{p+1}(-1)^i\varphi\big(\alpha(x_1),...,\hat{x_i},...,[x_i,x_j]_{\alpha},...,\alpha(x_{p+1}),\alpha(z)\big) \\
      &+& \sum_{i=1}^{p+1}(-1)^i\varphi\big(\alpha(x_1),...,\hat{x_i},...,\alpha(x_{p+1}),L(x_i)\cdot z\big)\nonumber
  \end{eqnarray}
  where $x_1,...,x_{p+1}\in \mathcal{L}(\mathcal{N}),\ z\in \mathcal{N}$ and $\hat{x}_i$ designed that $x_i$ is omitted.
  \end{definition}
  \begin{proposition}\label{propo4.4}
  If $\varphi\in C^p(\mathcal{N},\mathbb{K})$ be a $p$-cochain, then
  $$\delta^{p+1}\circ\delta^p(\varphi)=0$$
  \end{proposition}
  \begin{proof}Let  $\varphi$ be a $p$-cochain, $(x_i)_{1\leq i\leq p}\in \mathcal{L}(\mathcal{N})$ et $z\in \mathcal{N}$,
we can write $\delta^p$ and $\delta^{p+1}\circ\delta^p$ as
\begin{eqnarray*}
    & & \delta^p=\delta_1^p+\delta_2^p \\
  \textrm{and }& & \delta^{p+1}\circ\delta^p=\eta_{11}+\eta_{12}+\eta_{21}+\eta_{22}
\end{eqnarray*}  where $\eta_{ij}=\delta_i^{p+1}\circ\delta_j^p$, $1\leq i,j\leq 2$, and

\begin{eqnarray*}
  \delta_1^p\varphi(x_1,...,x_{p+1},z )&=& \sum_{1\leq i<j}^{p+1}(-1)^i\varphi\big(\alpha(x_1),...,\hat{x_i},...,[x_i,x_j]_{\alpha},...,\alpha(x_{p+1}),\alpha(z)\big) \\
  \delta_2^p\varphi(x_1,...,x_{p+1},z) &=& \sum_{i=1}^{p+1}(-1)^i\varphi\big(\alpha(x_1),...,\hat{x_i},...,\alpha(x_{p+1}),L(x_i)\cdot z\big)
\end{eqnarray*}

$\bullet$ Let us compute first $\eta_{11}\varphi(x_1,...,x_{p+1},z)$.  This is given by
\begin{eqnarray*}
    & & \eta_{11}(\varphi)(x_1,...,x_{p+1},z)\\
    &=& \sum_{1\leq i<k< j}^{p+1}(-1)^{i+k}\varphi\big(\alpha^2(x_1),...,\widehat{x_i},...,\widehat{\alpha(x_k)},...,[\alpha(x_k),[x_i,x_j]_\alpha]_\alpha,....,
\alpha^2(x_{p+1}),\alpha^2(z)\big) \\
    &+& \sum_{1\leq i<k< j}^{p+1}(-1)^{i+k-1}\varphi\big(\alpha^2(x_1),...,\widehat{\alpha(x_i)},...,\widehat{x_k},...,[\alpha(x_i),[x_k,x_j]_\alpha]_\alpha,....,
\alpha^2(x_{p+1}),\alpha^2(z)\big) \\
    &+& \sum_{1\leq i<k< j}^{p+1}(-1)^{i+k-1}\varphi\big(\alpha^2(x_1),...,\widehat{x_i},...,\widehat{[x_i,x_k]_\alpha},...,[[x_i,x_k]_\alpha,\alpha(x_j)]_\alpha,....,
\alpha^2(x_{p+1}),\alpha^2(z)\big).
\end{eqnarray*}
Whence applying the Hom-Leibniz identity \eqref{brackLei2} to $x_i,\ x_j,\ x_k\in \mathcal{L}(\mathcal{N})$, we find
$\eta_{11}=0$.\\
$\bullet$
\begin{eqnarray*}
    & & \eta_{21}(\varphi)(x_1,...,x_{p+1},z)+\eta_{12}(\varphi)(x_1,...,x_{p+1},z)=  \\
    & &  \sum_{1\leq i< j}^{p+1}(-1)^{i-1}\varphi\big(\alpha^2(x_1),...,\widehat{x}_i,...,\widehat{[x_i,x_j]_\alpha},...,
\alpha^2(x_{p+1}),L([x_i,x_j]_\alpha)\cdot\alpha(z)\big)
\end{eqnarray*}
and
\begin{eqnarray*}
    & & \eta_{22}(\varphi)(x_1,...,x_{p+1},z) \\
    &=& \sum_{1\leq i< j}^{p+1}(-1)^i\varphi\big(\alpha^2(x_1),...,\widehat{x}_i,...,\widehat{\alpha(x_j)},...,
\alpha^2(x_{p+1}),\big(L(\alpha(x_i))\cdot(L(x_j)\cdot z)\big)\big) \\
  &+& \sum_{1\leq i< j}^{p+1}(-1)^{i-1}\varphi\big(\alpha^2(x_1),...,\widehat{\alpha(x_i)},...,\widehat{x}_j,...,
\alpha^2(x_{p+1}),\big(L(\alpha(x_j))\cdot(L(x_i)\cdot z)\big)\big).
\end{eqnarray*}
Then applying  the Lemma \ref{3.1}  to $x_i,\ x_j\in \mathcal{L}(\mathcal{N})$ and $z\in \mathcal{N}$,
$\eta_{12}+\eta_{21}+\eta_{22}=0$.\\
Which ends the proof
  \end{proof}
  \begin{definition}The space of $p$-cocycles is defined by
  $$Z^p(\mathcal{N},\mathbb{K})=\{\varphi\in C^p(\mathcal{N},\mathbb{K}):\delta^p\varphi=0\}$$ and the space of $p$-coboundaries is defined by
  $$B^p(\mathcal{N},\mathbb{K})=\{\psi=\delta^{p-1}\varphi:\varphi\in C^{p-1}(\mathcal{N},\mathbb{K})\}$$
  \end{definition}
  \begin{lem}$B^p(\mathcal{N},\mathbb{K})\subset Z^p(\mathcal{N},\mathbb{K})$
  \end{lem}
  \begin{definition}We call $p^{\textrm{th}}$-cohomology group the quotient
  $$H^p(\mathcal{N},\mathbb{K})=\frac{Z^p(\mathcal{N},\mathbb{K})}{B^p(\mathcal{N},\mathbb{K})}$$
  \end{definition}
  \begin{example}
  Let $(\mathcal{N},[\cdot ,...,\cdot ])$ be a Nambu-Lie algebra $($see \cite{Fil} \cite{W.X}$)$ and $\{e_i\}_{i= 1}^{ n+1}$ be a basis such that
  \begin{equation}\label{filip}[e_1,...,\hat{e}_i,...,e_{n+1}]=(-1)^{i+1}\varepsilon_ie_i\ \ \ \textrm{or}\ \  [e_{i_1},....,e_{i_n}]=(-1)^n\sum_{i=1}^{n+1}\varepsilon_i\epsilon_{i_1,...,i_n}^ie_i\end{equation}
  where $\varepsilon_i=\pm1$ $($no sum over the $i$ of the $\varepsilon_i$ factors$)$ just introduce  signs that affect the different terms of the sum in $i$ and we have used  Filippov's notation.\\

  Note that we might equally well have the $\epsilon_{i_1,...,i_n}^i$ without signs $\varepsilon_i$ in $\ref{filip}$ by taking $\epsilon_{i_1,...,i_n}^i=\eta^{ij}\epsilon_{i_1,...,i_n,j}$, where $\epsilon_{1,...,n,(n+1)}=1$ and $\eta$ is a $(n+1)\times(n+1)$ diagonal matrix with $+1$ and $-1$ in  places indicated by the $\varepsilon_i$'s. We shall keep nevertheless the customary $\varepsilon_i$ factors above as  in $e.g.$ \cite{W.X}.\\

  Let $\alpha :\mathcal{N}\rightarrow \mathcal{N}$ be a morphism of Nambu-Lie algebras. Then using Theorem 1.5, $\mathcal{N}_\alpha=(\mathcal{N},[\cdot ,...,\cdot ]_\alpha,\tilde{\alpha}=(\alpha,...,\alpha))$ is a Hom-Nambu-Lie algebra where the bracket $[\cdot ,...,\cdot ]_\alpha$ is given by
  \begin{equation}[e_1,...,\hat{e}_i,...,e_{n+1}]_\alpha=(-1)^{i+1}\varepsilon_i\alpha(e_i)\ \ \ \textrm{or}\ \  [e_{i_1},....,e_{i_n}]_\alpha=(-1)^n\sum_{i=1}^{n+1}\varepsilon_i\epsilon_{i_1,...,i_n}^i\alpha(e_i).\end{equation}

We establish the following result.
    \begin{lem}Any $1$-cochain of the Hom-Nambu-Lie algebra $\mathcal{N}_\alpha$ is a $1$-coboundary (and thus a trivial $1$-cocycle).
    \end{lem}
    \begin{proof}
    Let $\varphi\in C^1(\mathcal{N},\mathbb{K})$ be a $1$-cochain on $\mathcal{N}_\alpha$, $\varphi$ is determined by its coordinates $\varphi_{i_1,...,i_n}=\varphi(e_{i_1},...,e_{i_n})$. We now show that, in fact, a $1$-cochain on $\mathcal{N}_\alpha$ is a $1$-coboudary, that is  there exists a $0$-cochain $\phi$ such that
    \begin{equation}\label{7.4}
    \varphi_{i_1,...,i_n}=-\phi([e_{i_1},...,e_{i_n}])=-\sum_{k=1}^{n+1}\varepsilon_k\epsilon_{i_1,...,i_n}^k\phi_k,
    \end{equation}
    where $\phi_k=\phi\circ\alpha(e_k)$. Indeed, given $\varphi$ then the $0$-cochain $\phi$ is given by
    \begin{equation}
    \phi_k=-\frac{\varepsilon_k}{n!}\sum_{i_1...i_n}^{n+1}\epsilon_k^{i_1,...,i_n}\varphi_{i_1,...,i_n}
    \end{equation}
    has the desired property \eqref{7.4}:
    \begin{eqnarray}
      -\phi([e_{i_1},...,e_{i_n}])&=&-\sum_{k=1}^{n+1}\varepsilon_k\epsilon_{i_1,...,i_n}^k\phi_k\nonumber \\
        &=& \sum_{k=1}^{n+1}\epsilon_{i_1,...,i_n}^k \frac{\varepsilon_k^2}{n!}\sum_{j_1...j_n}^{n+1}\epsilon_k^{j_1,...,j_n}\varphi_{j_1,...,j_n}\nonumber\\
        &=& \frac{1}{n!}\sum_{j_1...j_n}^{n+1}\epsilon_{i_1,...,i_n}^{j_1,...,j_n}\varphi_{j_1,...,j_n}=\varphi_{i_1,...,i_n}
    \end{eqnarray}
    which proves the lemma.
    \end{proof}
  \end{example}
  \section{Deformation of $n$-ary Hom-Nambu-Lie algebras }

  Let $\mathbb{K}[[t]]$ be the power series ring in one variable $t$ and coefficients in $\mathbb{K}$ and $\mathcal{N}[[t]]$ be the set of formal series whose coefficients are elements of the vector space $\mathcal{N}$, ($\mathcal{N}[[t]]$ is obtained by extending the coefficients domain of $\mathcal{N}$ from $\mathbb{K}$ to $\mathbb{K}[[t]]$). Given a $\mathbb{K}$-$n$-linear map $\varphi:\mathcal{N}\times...\times \mathcal{N}\rightarrow \mathcal{N}$, it admits naturally an extension to a  $\mathbb{K}[[t]]$-$n$-linear map $\varphi:\mathcal{N}[[t]]\times...\times \mathcal{N}[[t]]\rightarrow \mathcal{N}[[t]]$, that is, if $x_i=\dl\sum_{j\geq0}a_i^jt^j$, $1\leq i\leq n$ then $\varphi(x_1,...,x_n)=\dl\sum_{j_1,...,j_n\geq0}t^{j_1+...+j_n}\varphi(a_1^{j_1},...,a_n^{j_n})$. The same holds for linear map.

  \begin{definition}Let $(\mathcal{N},[\cdot ,...,\cdot ],\widetilde{\alpha}),\ \widetilde{\alpha}=(\alpha_i)_{1\leq i\leq n-1}$ be a Hom-Nambu-Lie algebra. A formal deformation of the Hom-Nambu-Lie algebra $\mathcal{N}$ is given by a $\mathbb{K}[[t]]$-$n$-linear map $$[\cdot ,...,\cdot ]_t:\mathcal{N}[[t]]\times...\times \mathcal{N}[[t]]\rightarrow \mathcal{N}[[t]]$$of the form $[\cdot ,...,\cdot ]_t=\dl\sum_{i\geq0}t^i[\cdot ,...,\cdot ]_i$ where each $[\cdot ,...,\cdot ]_i$ is a $\mathbb{K}[[t]]$-$n$-linear map $[\cdot ,...,\cdot ]_i:\mathcal{N}\times...\times \mathcal{N}\rightarrow \mathcal{N}$ (extending to be  $\mathbb{K}[[t]]$-$n$-linear), and $[\cdot ,...,\cdot ]_0=[\cdot ,...,\cdot ]$ such that for $(x_i)_{1\leq i\leq n-1},\ (y_i)_{1\leq i\leq n}\in \mathcal{N}$
  $$
  \big[\alpha_1(x_1),....,\alpha_{n-1}(x_{n-1}),[y_1,....,y_n]_t\big]_t= $$
  \begin{equation}\sum_{i=1}^{n-1}\big[\alpha_1(y_1),....,\alpha_{i-1}(y_{i-1}),[x_1,....,x_{n-1},y_i]_t
  ,\alpha_i(y_{i+1}),...,\alpha_{n-1}(y_n)\big]_t.
  \end{equation}
  The deformation is said to be of order $k$ if $[\cdot ,...,\cdot ]_t=\dl\sum_{i=0}^kt^i[\cdot ,...,\cdot ]_i$ and infinitesimal if $t^2=0$.\\
  In terms of elements $x=(x_i)_{1\leq i\leq n-1},\ y=(y_i)_{1\leq i\leq n-1}\in \mathcal{L}(\mathcal{N})$ and setting $z=y_n$ the above condition reads
  \begin{equation}\label{8.1} L_t([x ,y ]_{\alpha})\cdot \alpha_n(z)=L_t(\tilde{\alpha}(x))\cdot \big(L_t(y)\cdot z\big)-L_t(\tilde{\alpha}(y))\cdot \big(L_t(x).z\big)\end{equation}
  where $L_t(x)\cdot z=[x_1,...,x_{n-1},z]_t$ and $\tilde{\alpha}(x)=(\alpha_i(x_i))_{1\leq i\leq n-1}$.\\
  Now let $(\mathcal{N},[\cdot ,...,\cdot ],\alpha)$ be a multiplicative Hom-Nambu-Lie (i.e. $\alpha_1=...=\alpha_n=\alpha$).\\
   Eq. \eqref{8.1} implies, keeping only terms linear in $t$,
  \begin{eqnarray*}
      & &\big[\alpha(x_1),....,\alpha(x_{n-1}),\psi(y_1,....,y_n)\big] + \psi\big(\alpha(x_1),....,\alpha(x_{n-1}),[y_1,....,y_n]\big) \\
      &=&  \sum_{i=1}^n \big[\alpha(y_1),....,\alpha(y_{i-1}),\psi(x_1,....,x_{n-1},y_i)
  ,\alpha(y_{i+1}),...,\alpha(y_n)\big]  \\
      &+&\sum_{i=1}^n\psi\big(\alpha(y_1),....,\alpha(y_{i-1}),[x_1,....,x_{n-1},y_i]
  ,\alpha(y_{i+1}),...,\alpha(y_n)\big).
  \end{eqnarray*}
This expression may be read as the $1$-cocycle condition $\delta^1\psi=0$ for the $\mathcal{N}$-valued cochain $\psi$. In terms of $x,\ y\in \mathcal{L}(\mathcal{N})$ it may be written, (setting again $y_n=z$ ), as
\begin{eqnarray}
    \delta^1\psi(x,y,z) &=& \psi(\alpha(x),L(y)\cdot z)-
  \psi(\alpha(y),L(x)\cdot z)-\psi([x_1,x_2]_{\alpha},\alpha(z)) \\
    &+&L(\alpha(x))\cdot \psi(y,z)- L(\alpha(y))\cdot \psi(x,z)+\big(\psi(x,\ \ )\cdot y\big)\bullet_{\alpha}\alpha(z)\nonumber
  \end{eqnarray}
  where \begin{equation}\big(\psi(x,\ \ )\cdot y\big)\bullet_{\alpha}\alpha(z)=\dl\sum_{i=0}^{n-1}[\alpha(y_1),...,\psi(x,y_i),...,\alpha(y_{n-1}),\alpha(z)].\end{equation}
  \end{definition}
  \begin{definition}
a $p$-cochains is an $p+1$-linear map
 $
  \varphi:\mathcal{L}(\mathcal{N})\otimes...\otimes \mathcal{L}(\mathcal{N})\wedge \mathcal{N} \longrightarrow \mathcal{N} $,
such that $$\alpha\circ\varphi(x_1,...,x_p,z)=\varphi(\alpha(x_1),...,\alpha(x_p),\alpha(z)).$$
We denote the set of a $p$-cochain by $C^p(\mathcal{N},\mathcal{N})$
\end{definition}
  \begin{definition}We call, for $p\geq 1$, $p$-coboundary operator of the multiplicative Hom-Nambu-Lie $(\mathcal{N},[\cdot ,...,\cdot ],\alpha)$ the linear map $\delta^p:C^p(\mathcal{N},\mathcal{N})\rightarrow C^{p+1}(\mathcal{N},\mathcal{N})$ defined by

  \begin{eqnarray}
                                  \delta^p\psi(x_1,...,x_p,x_{p+1},z) &=&\sum_{1\leq i\leq j}^{p+1}(-1)^i\psi\big(\alpha(x_1),...,\widehat{\alpha(x_i)},...,
                                  \alpha(x_{j-1}),[x_i,x_j]_{\alpha},...,\alpha(x_{p+1}),\alpha(z)\big)\nonumber\\
                                    &+& \sum_{i=1}^{p+1}(-1)^i\psi\big(\alpha(x_1),...,\widehat{\alpha(x_i)},...,
                                    \alpha(x_{p+1}),L(x_i)\cdot z\big) \nonumber\\
                                    &+& \sum_{i=1}^{p+1}(-1)^{i+1}L(\alpha^p(x_i))\cdot \psi\big(x_1,...,\widehat{x_i},...,
                                  x_{p+1},z\big)\nonumber\\
                                  &+&  (-1)^p\big(\psi(x_1,...,x_p,\ )\cdot x_{p+1}\big)\bullet_{\alpha}\alpha^p(z)
                                \end{eqnarray}
where
  \begin{equation}\big(\psi(x_1,...,x_p,\ )\cdot x_{p+1}\big)\bullet_{\alpha}\alpha^p(z)=\dl\sum_{i=1}^{n-1}[\alpha^p(x_{p+1}^1),...,\psi(x_1,...,x_p,x_{p+1}^i ),...,\alpha^p(x_{p+1}^{n-1}),\alpha^p(z)],\end{equation}\\for all $x_i=(x_i^j)_{1\leq j\leq n-1}\in \mathcal{L}(\mathcal{N}),\ 1\leq i\leq p+1$, $z\in \mathcal{N}$ and $\widehat{x}_i$ designed that $x_i$ is omitted.

  \end{definition}
  \begin{proposition}\label{opercob}Let $\psi\in C^p(\mathcal{N},\mathcal{N})$ be a $p$-cochain then
  $$\delta^{p+1}\circ\delta^p(\psi)=0.$$
  \end{proposition}
  \begin{proof}Let $\psi$ be a $p$-cochain, $x_i=(x_i^j)_{1\leq j\leq n-1}\in \mathcal{L}(\mathcal{N}),\ 1\leq i\leq p+2$ and $z\in \mathcal{N}$  we can write $\delta^p$ and $\delta^{p+1}\circ\delta^p$ as
  \begin{eqnarray*}
      & & \delta^p=\delta_1^p+\delta_2^p+\delta_3^p+\delta_4^p, \\
    \textrm{and} & & \delta^{p+1}\circ\delta^p= \sum_{i,j=1}^4\eta_{ij},
  \end{eqnarray*}
  when $\eta_{ij}=\delta_i^{p+1}\circ\delta_j^p$ and
  \begin{eqnarray*}
    \delta_1^p\psi(x_1,...,x_{p+1},z )&=& \sum_{1\leq i<j}^{p+1}(-1)^i\psi\big(\alpha(x_1),...,\widehat{x_i},...,[x_i,x_j]_{\alpha},...,\alpha(x_{p+1}),\alpha(z)\big) \\
  \delta_2^p\psi(x_1,...,x_{p+1},z) &=& \sum_{i=1}^{p+1}(-1)^i\psi\big(\alpha(x_1),...,\widehat{x_i},...,\alpha(x_{p+1}),L(x_i).z\big) \\
    \delta_3^p\psi(x_1,...,x_{p+1},z) &=& \sum_{i=1}^{p+1}(-1)^{i+1}L(\alpha^p(x_i))\cdot \psi\big(x_1,...,\widehat{x_i},...,
                                  x_{p+1},z\big) \\
    \delta_4^p\psi(x_1,...,x_{p+1},z) &=&  (-1)^p\big(\psi(x_1,...,x_p,\ )\cdot x_{p+1}\big)\bullet_{\alpha}\alpha^p(z)
  \end{eqnarray*}
To simplify the notations we replace $L(x)\cdot z$ by $x\cdot z$.\\
  The proof that $\eta_{11}+\eta_{12}+\eta_{21}+\eta_{22}=0$ is similar to the proof in Proposition \ref{propo4.4}.\\
  On the other hand, we have
\begin{eqnarray*}
            \star \ \eta_{13}\psi(x_1,...,x_{p+2},z )
              &=&\sum_{1\leq i<j<k}^{p+2}  \big\{(-1)^{k+i}\alpha^{p+1}(x_k)\cdot \psi(\alpha(x_1),...,\widehat{x}_i,...,[x_i,x_j]_\alpha,...,\widehat{x}_k,...,\alpha(z)) \\
                &+& (-1)^{j+i}\alpha^{p+1}(x_j)\cdot \psi(\alpha(x_1),...,\widehat{x}_i,...,\widehat{x}_j,...,[x_i,x_k]_\alpha,...,\alpha(z)) \\
                &+& (-1)^{j+i-1}\alpha^{p+1}(x_i)\cdot \psi(\alpha(x_1),...,\widehat{x}_i,...,\widehat{x}_j,...,[x_j,x_k]_\alpha,...,\alpha(z))\big\}\\
            & & \\
            \star\ \eta_{31}\psi(x_1,...,x_{p+2},z ) &=& -\eta_{13}\psi(x_1,...,x_{p+2},z ) \\
                & +& \sum_{1\leq i<j}^{p+2}  (-1)^{i+j}\alpha^p([x_i,x_j]_\alpha)\cdot \alpha\big(\psi(x_1,...,\widehat{x}_i,...,\widehat{x}_j,...,z)\big)\\
                & & \\
                \star\ \eta_{33}\psi(x_1,...,x_{p+2},z )
              &=&\sum_{1\leq i<j}^{p+2} \Big\{  (-1)^{i+j}\alpha^{p+1}(x_i)\cdot \big(\alpha^p(x_j).\big(\psi(x_1,...,\widehat{x}_i,...,\widehat{x}_j,...,z)\big)\big) \\
                &+& (-1)^{i+j-1}\alpha^{p+1}(x_j)\cdot \big(\alpha^p(x_i)\cdot \big(\psi(x_1,...,\widehat{x}_i,...,\widehat{x}_j,...,z)\big)\big)\Big\}
            \end{eqnarray*}
Then, applying  Lemma \ref{3.1} to  $\alpha^p(x_i)\in \mathcal{L}(\mathcal{N})$, $\alpha^p(x_j)\in \mathcal{L}(\mathcal{N})$ et $\psi(x_1,...,\widehat{x}_i,...,\widehat{x}_j,...,z)\in \mathcal{N}$, we have$$\eta_{13}+\eta_{33}+\eta_{31}=0.$$
by the same calculation, we can prove that
$$\eta_{23}+\eta_{32}=0.$$
  \begin{eqnarray*}
    & &\star\ \eta_{14}\psi(x_1,...,x_{p+2},z) \\&=& (-1)^p\sum_{1\leq i<j}^{p+1}(-1)^i\sum_{k=1}^{n-1}\big[\alpha^{p+1}(x_{p+2}^1),...,    \psi\big(\alpha(x_1),...,\widehat{x}_i,...,[x_i,x_j]_\alpha,...,\alpha(x_p),\alpha(x_{p+1}^k)\big),...,\alpha^{p+1}(z)\big] \\
      &+& (-1)^p\sum_{i=1}^{p+1}(-1)^i\sum_{k,l=1;k\neq l}^{n-1}\big[\alpha^{p+1}(x_{p+2}^1),...,\alpha^p(x_i\cdot x_{p+2}^l),    ...,\psi\big(\alpha(x_1),...,\widehat{x}_i,...,\alpha(x_{p+1}),\alpha(x_{p+2}^k)\big),...,\alpha^{p+1}(z)\big] \\
      &+& (-1)^p\sum_{i=1}^{p+1}(-1)^i\sum_{k=1}^{n-1}\big[\alpha^{p+1}(x_{p+2}^1),...,    \psi\big(\alpha(x_1),...,\widehat{x}_i,...,\alpha(x_{p+1}),x_i\cdot x_{p+2}^k\big),...,\alpha^{p+1}(z)\big].
      \end{eqnarray*}
The first term in $\eta_{14}$  is equal to $-\eta_{41}$, hence
      \begin{eqnarray*}
  & & \star\ (\eta_{14}+\eta_{41})\psi(x_1,...,x_{p+2},z)\\
      &=& \sum_{i=1}^{p+1}(-1)^{p+i}\sum_{k,l=1;k\neq l}^{n-1}\big[\alpha^{p+1}(x_{p+2}^1),...,\alpha^p(x_i\cdot x_{p+2}^l),...,    \psi\big(\alpha(x_1),...,\widehat{x}_i,...,\alpha(x_{p+1}),\alpha(x_{p+2}^k)\big),...,\alpha^{p+1}(z)\big] \\
      &+& \sum_{i=1}^{p+1}(-1)^{p+i}\sum_{k=1}^{n-1}\big[\alpha^{p+1}(x_{p+2}^1),...,    \psi\big(\alpha(x_1),...,\widehat{x}_i,...,\alpha(x_{p+1}),x_i\cdot x_{p+2}^k\big),...,\alpha^{p+1}(z)\big]\\
    \textrm{and}& & \\
  & & \star\ \eta_{24}\psi(x_1,...,x_{p+2},z)\\ &=& \sum_{i=1}^{p+1}\sum_{k=1}^{n-1}(-1)^{p+i}\big[\alpha^{p+1}(x_{p+2}^1),...,    \psi\big(\alpha(x_1),...,\widehat{x}_i,...,\alpha(x_{p+1}),\alpha(x_{p+2}^k)\big),...,\alpha^p(x_i\cdot z)\big]\\
    &+& \sum_{k=1}^{n-1}\big[\alpha^{p+1}(x_{p+1}^1),...,    \psi\big(\alpha(x_1),...,\widehat{x}_i,...,\alpha(x_p),\alpha(x_{p+1}^k)\big),...,\alpha^p(x_{p+2}\cdot z)\big]\\
    & & \\
    & &\star\ \eta_{42}\psi(x_1,...,x_{p+2},z)\\
      &=& (-1)^{p+1}\sum_{i=1}^{p+1}(-1)^i\sum_{k=1}^{n-1}\big[\alpha^{p+1}(x_{p+2}^1),...,  \psi\big(\alpha(x_1),...,\widehat{x}_i,...,\alpha(x_{p+1}),x_i\cdot x_{p+2}^k\big),...,\alpha^{p+1}(z)\big].
      \end{eqnarray*}
  Hence,  $-\eta_{42}$ and the second term of $ (\eta_{14}+\eta_{41})$ are equal.\\

    Using the Hom-Nambu identity for any integers $1\leq i\leq p+1$ et $1\leq k\leq n-1$
\begin{eqnarray*}
    & & \alpha^{p+1}(x_i)\cdot \big[\alpha^p(x_{p+2}^1),...,    \psi\big(x_1,...,\widehat{x}_i,...,x_{p+1},x_{p+2}^k\big),...,\alpha^p(z)\big] \\
    &= & \sum_{l=1;l\neq k}^{n-1}\Big\{\big[\alpha^{p+1}(x_{p+2}^1),...,\alpha^p(x_i\cdot x_{p+2}^l),...,
        \psi\big(\alpha(x_1),...,\widehat{x}_i,...,\alpha(x_{p+1}),\alpha(x_{p+2}^k)\big),...,\alpha^{p+1}(z)\big]\Big\} \\
    & &+\big[\alpha^{p+1}(x_{p+2}^1),...,
      \psi\big(\alpha(x_1),...,\widehat{x}_i,...,\alpha(x_{p+1}),\alpha(x_{p+2}^k)\big),...,\alpha^p(x_i\cdot z)\big]  \\
    & & +\big[\alpha^{p+1}(x_{p+1}^1),...,\alpha^p(x_i)\cdot  \psi\big(x_1,...,\widehat{x}_i,...,x_{p+1},x_{p+2}^k\big),...,\alpha^{p+1}(z)\big]
\end{eqnarray*}
    when we add the four terms $\eta_{14}$, $\eta_{41}$, $\eta_{24}$ and $\eta_{42}$, we have the following  expression
    \begin{eqnarray*}
    & & (\eta_{14}+\eta_{41}+\eta_{24}+\eta_{42})\psi(x_1,...,x_{p+2},z)\\
      &=& \sum_{i=1}^{p+1}(-1)^{i+p}\sum_{l=1}^{n-1}\big[\alpha^{p+1}(x_{p+2}^1),...,  \alpha^p(x_i)\cdot  \psi\big(\alpha(x_1),...,\widehat{x}_i,...,\alpha(x_{p+1}),\alpha(x_{p+2}^k)\big),...,\alpha^{p+1}(z)\big] \\
      &+& (-1)^{p-1}\sum_{i=1}^{p+1}(-1)^i\sum_{k=1}^{n-1}\alpha^{p+1}(x_i)\cdot \big[\alpha^p(x_{p+2}^1),...,    \psi\big(x_1,...,\widehat{x}_i,...,x_{p+1},x_{p+2}^k\big),...,\alpha^p(z)\big]\\
    &+& \sum_{k=1}^{n-1}\big[\alpha^{p+1}(x_{p+1}^1),...,    \psi\big(\alpha(x_1),...,\widehat{x}_i,...,\alpha(x_p),\alpha(x_{p+1}^k)\big),...,\alpha^p(x_{p+2}\cdot z)\big]\\
      \textrm{and}& & \\
      & &\star\ \eta_{43}\psi(x_1,...,x_{p+2},z)\\
            &=& \sum_{i=1}^{p+1}(-1)^{p+i}\sum_{k=1}^{n-1}\alpha^{p+1}(x_i)\cdot \big[\alpha^{p}(x_{p+2}^1),...,    \psi\big(x_1,...,\widehat{x}_i,...,x_{p+1},x_{p+2}^k\big),...,\alpha^{p}(z)\big] \\
  &-& \sum_{k=1}^{n-1}\alpha^{p+1}(x_{p+2})\cdot \big[\alpha^p(x_{p+1}^1),...,  \psi\big(x_1,...,x_p,x_{p+1}^k\big),...,\alpha^p(z)\big],
  \end{eqnarray*}
      \begin{eqnarray*}
  & & \eta_{34}\psi(x_1,...,x_{p+2},z)\\
  &=& \sum_{i=1}^{p+1}(-1)^{i+p+1}\sum_{l=1}^{n-1}\big[\alpha^{p+1}(x_{p+2}^1),...,  \alpha^p(x_i)\cdot  \psi\big(\alpha(x_1),...,\widehat{x}_i,...,\alpha(x_{p+1}),\alpha(x_{p+2}^k)\big),...,\alpha^{p+1}(z)\big].
      \end{eqnarray*}
      Hence
      \begin{eqnarray*}
  & &  (\eta_{14}+\eta_{41}+\eta_{24}+\eta_{42}+\eta_{34}+\eta_{43})\psi(x_1,...,x_{p+2},z)=-\eta_{44}\psi(x_1,...,x_{p+2},z)\\
            &=& -\sum_{i=1}^{n-1}\sum_{k=1}^{n-1}\big[\alpha^{p+1}(x_{p+2}^1),...,  [\alpha^p(x_{p+1}^1),..., \psi\big(x_1,...,x_p,x_{p+1}^k\big),...,\alpha^p(x_{p+2}^i)],...,\alpha^{p+1}(x_{p+2}^{n-1}),\alpha^{p+1}(z)\big]\\
            &=&-\sum_{k=1}^{n-1}\alpha^{p+1}(x_{p+2})\cdot [\alpha^p(x_{p+1}^1),..., \psi\big(x_1,...,x_p,x_{p+1}^k\big),...,\alpha^p(z)]\\
            &+&\sum_{k=1}^{n-1}[\alpha^{p+1}(x_{p+1}^1),..., \psi\big(\alpha(x_1),...,\alpha(x_p),\alpha(x_{p+1}^k)\big),...,\alpha^p(x_{p+2}\cdot z)].
\end{eqnarray*}
Then, we have $$\eta_{14}+\eta_{41}+\eta_{24}+\eta_{42}+\eta_{34}+\eta_{43}+\eta_{44}=0,$$
which ends the proof.
  \end{proof}
  \begin{definition}The space of $p$-cocycles is defined by
  $$Z^p(\mathcal{N},\mathcal{N})=\{\varphi\in C^p(\mathcal{N},\mathcal{N}):\delta^p\varphi=0\},$$ and the space of $p$-coboundaries is defined by
  $$B^p(\mathcal{N},\mathcal{N})=\{\psi=\delta^{p-1}\varphi:\varphi\in C^{p-1}(\mathcal{N},\mathcal{N})\}.$$
  \end{definition}
  \begin{lem}$B^p(\mathcal{N},\mathcal{N})\subset Z^p(\mathcal{N},\mathcal{N})$
  \end{lem}
  \begin{definition}We call the $p^{\textrm{th}}$-cohomology group  the quotient
  $$H^p(\mathcal{N},\mathcal{N})=\frac{Z^p(\mathcal{N},\mathcal{N})}{B^p(\mathcal{N},\mathcal{N})}.$$
  \end{definition}
\section{Cohomology of $n$-ary Hom-algebras induced by cohomology of Hom-Leibniz algebras}
\subsection{Cohomology of ternary Hom-Nambu algebras induced by cohomology of Hom-Leibniz algebras}
In this section we extend to ternary multiplicative Hom-Nambu-Lie algebras the Takhtajan's construction of a cohomology of ternary Nambu-Lie algebras starting from Chevalley-Eilenberg cohomology of binary Lie algebras, (see \cite{Takhtajan0,Takhtajan1,Takhtajan2}). The cohomology of multiplicative Hom-Lie algebras was introduced in \cite{AEM} and independently in \cite{Sheng}.

The cohomology complex for Leibniz algebras was defined by Loday-Pirashvili in \cite{lodayPirashvili93}. We extend it to Hom-Leibniz algebras as follows.

Let $(A, [\cdot , \cdot ], \alpha )$ be a Hom-Leibniz algebras and $\mathcal{C}_\mathcal{L}(A,A)$ be the set of cochains
$\mathcal{C}^p_\mathcal{L}(A,A)=Hom(\otimes^p A,A)$ for $n\geq 1$. We set $\mathcal{C}^0_\mathcal{L}(A,A)=A$. We define a coboundary operator $d$ by $d\varphi(a)=-[\varphi, a]$ when $\varphi\in \mathcal{C}^0_\mathcal{L}(A,A)$ and for $p\geq 1$, $\varphi\in \mathcal{C}^p_\mathcal{L}(A,A)$, $a_1,\cdots ,a_{p+1}\in A$

\begin{align}\label{Leibnizcohomo}
d^p\varphi(a_1,\cdots ,a_{p+1})& =\sum_{k=1}^p(-1)^{k-1} \big[ \alpha^{p-1}(a_k),\varphi(a_1,\cdots,\widehat{a_k},\cdots,a_{p+1})\big]\\
\ & +(-1)^{p+1} \big[ \varphi(a_1\otimes\cdots\otimes a_p),\alpha^{p-1}(a_{p+1}) \big]  \nonumber\\
\ & +\sum_{1\leq k<j}^{p+1}{(-1)^k \varphi(\alpha(a_1)\otimes \cdots \otimes\widehat{a_k}\otimes \cdots\otimes\alpha(a_{j-1})\otimes [a_k,a_j]\otimes\alpha(a_{j+1})\otimes \cdots\otimes\alpha(a_{p+1}))}\nonumber
\end{align}
Notice that we recover the classical case when $\alpha =id.$\\

We aim now to derive the cohomology of a ternary Hom-Nambu algebra from the cohomology of Hom-leibniz algebra following the procedure described for ternary Nambu algebra in \cite{Takhtajan0}.\\

Let $(\mathcal{N},[\cdot ,\cdot,\cdot ],\alpha)$ be a multiplicative ternary Hom-Nambu-Lie algebra. Using Proposition \ref{HomLeibOfHomNambu} the triple $(\mathcal{L}(\mathcal{N})=\mathcal{N}\otimes\mathcal{N},[\cdot,\cdot]_\alpha,\alpha)$  where the bracket is defined for $x=x_1\otimes x_2$ and $y=y_1\otimes y_2$ by
\begin{equation}[x,y]=[x_1,x_2,y_1]\otimes \alpha (y_2)+\alpha (y_1)\otimes [x_1,x_2,y_2],
\end{equation}
is a  Hom-Leibniz  algebra.

\begin{thm}Let $(\mathcal{N},[\cdot ,\cdot,\cdot ],\alpha)$ be a multiplicative ternary Hom-Nambu-Lie algebra and $\mathcal{C}^p_\mathcal{N}(\mathcal{N},\mathcal{N})=Hom(\otimes^{2p+1}\mathcal{N},\mathcal{N})$ for $n\geq 1$ be the cochains. Let
$\Delta :\mathcal{C}^p_\mathcal{N}(\mathcal{N},\mathcal{N})\rightarrow \mathcal{C}^{p+1}_\mathcal{L}(\mathcal{L},\mathcal{L})$ be the linear map defined  for $p=0$ by $$\Delta\varphi(x_1\otimes x_2)=x_1\otimes\varphi( x_2)+\varphi(x_1)\otimes x_2$$ and for $p>0$
\begin{eqnarray*}(\Delta\varphi )(a_1,\cdots ,a_{p+1})=\alpha^{p-1}(x_{2p+1})\otimes\varphi(a_1,\cdots ,a_{p}\otimes x_{2p+2})+\varphi(a_1,\cdots ,a_{p}\otimes x_{2p+1})\otimes \alpha^{p-1}(x_{2p+2}),
\end{eqnarray*}
where we set $a_j=x_{2j-1}\otimes x_{2j}.$

Then there exists a cohomology complex $(\mathcal{C}^\bullet_\mathcal{N}(\mathcal{N},\mathcal{N}),\delta )$ for ternary Hom-Nambu-Lie algebras such that $$d\circ \Delta =\Delta\circ \delta.$$

The coboundary map $\delta: \mathcal{C}^p_\mathcal{N}(\mathcal{N},\mathcal{N})\rightarrow \mathcal{C}^{p+1}_\mathcal{N}(\mathcal{N},\mathcal{N})$ is defined for $\varphi\in \mathcal{C}^p_\mathcal{N}(\mathcal{N},\mathcal{N})$ by
\begin{align}
\delta\varphi (x_1\otimes\cdots \otimes x_{2p+1})& =\sum_{j=1}^p{\sum_{k=2j+1}^{2p+1}{(-1)^j\varphi(\alpha(x_1)\otimes\cdots \otimes [x_{2j-1},x_{2j},x_k]\otimes \cdots \otimes \alpha(x_{2p+1}))}}+\\
\ & \sum_{k=1}^{p}{(-1)^{k-1}[\alpha^{p-1}(x_{2k-1}),\alpha^{p-1}(x_{2k}),\varphi(x_1\otimes\cdots \otimes\widehat{x_{2k-1}}\otimes \widehat{x_{2k}}\otimes \cdots \otimes x_{2p+1})]}+
\nonumber\\
\ &(-1)^{n+1}[\alpha^{p-1}(x_{2p-1}),\varphi(x_1\otimes  \cdots \otimes x_{2p-2}\otimes x_{2p}),\alpha^{p-1}(x_{2p+1})]+
\nonumber\\
\nonumber \ &(-1)^{p+1}[\varphi(x_1\otimes \cdots \otimes x_{2p-1} ),\alpha^{p-1}(x_{2p}),\alpha^{p-1}(x_{2p+1})]
\end{align}
\end{thm}
\begin{proof}The proof is a particular case of Theorem \ref{constakh} proof.
\end{proof}
\begin{remark}
The theorem shows that one may derive the cohomology complex of ternary Hom-Nambu-Lie algebras from the cohomology complex of Hom-Leibniz algebras.
\end{remark}

\subsection{Cohomology of n-ary  Hom-Nambu-Lie algebras induced by cohomology of Hom-Leibniz algebras}
We generalize in this section the result of the previous section to $n$-ary Hom-Nambu-Lie  algebras.\\

Let $(\mathcal{N},[\cdot ,...,\cdot ],\alpha)$ be a multiplicative $n$-ary Hom-Nambu-Lie algebra and   the triple $(\mathcal{L}(\mathcal{N})=\mathcal{N}^{\otimes n-1},[\cdot,\cdot]_\alpha,\alpha)$  be the Hom-Leibniz algebra associates to $\mathcal{N}$ where the bracket is defined in \eqref{brackLei}.
\begin{thm}\label{constakh}Let $(\mathcal{N},[\cdot ,...,\cdot ],\alpha)$ be a multiplicative $n$-ary Hom-Nambu-Lie algebra and $\mathcal{C}^p_\mathcal{N}(\mathcal{N},\mathcal{N})=Hom(\otimes^{p}\mathcal{L}(\mathcal{N})\otimes\mathcal{N},\mathcal{N})$ for $p\geq 1$ be the cochains. Let
$\Delta :\mathcal{C}^p_\mathcal{N}(\mathcal{N},\mathcal{N})\rightarrow \mathcal{C}^{p+1}_\mathcal{L}(\mathcal{L},\mathcal{L})$ be the linear map defined for $p=0$ by \begin{eqnarray}\Delta\varphi(x_1\otimes\cdots\otimes x_{n-1})=\dl\sum_{i=0}^{n-1}x_1\otimes\cdots\otimes\varphi(x_i)\otimes\cdots\otimes x_{n-1}\end{eqnarray} and for $p>0$ by
\begin{eqnarray}\label{defdelta}(\Delta\varphi )(a_1,\cdots ,a_{p+1})=\sum_{i=1}^{n-1}\alpha^{p-1}(x_{p+1}^1)\otimes\cdots\otimes\varphi(a_1,\cdots ,a_{n-1}\otimes x_{p+1}^i)\otimes\cdots\otimes \alpha^{n-1}(x_{p+1}^{n-1}),
\end{eqnarray}
where we set $a_j=x_{j}^1\otimes\cdots\otimes x_j^{n-1}.$

Then there exists a cohomology complex $(\mathcal{C}^\bullet_\mathcal{N}(\mathcal{N},\mathcal{N}),\delta )$ for $n$-ary Hom-Nambu-Lie algebras such that $$d\circ \Delta =\Delta\circ \delta.$$

The coboundary map $\delta: \mathcal{C}^p_\mathcal{N}(\mathcal{N},\mathcal{N})\rightarrow \mathcal{C}^{p+1}_\mathcal{N}(\mathcal{N},\mathcal{N})$ is defined for $\varphi\in \mathcal{C}^p_\mathcal{N}(\mathcal{N},\mathcal{N})$ by
\begin{eqnarray*}
                                  \delta\varphi(a_1,...,a_p,a_{p+1},x) &=&\sum_{1\leq i\leq j}^{p+1}(-1)^i\varphi\big(\alpha(a_1),...,\widehat{\alpha(a_i)},...,
                                  \alpha(a_{j-1}),[a_i,a_j]_{\alpha},...,\alpha(a_{p+1}),\alpha(x)\big)\\
                                    &+& \sum_{i=1}^{p+1}(-1)^i\varphi\big(\alpha(a_1),...,\widehat{\alpha(a_i)},...,
                                    \alpha(a_{p+1}),L(a_i).x\big) \\
                                    &+& \sum_{i=1}^{p+1}(-1)^{i+1}L(\alpha^p(a_i))\cdot \varphi\big(a_1,...,\widehat{a_i},...,
                                  a_{p+1},x\big)\\
                                  &+&  (-1)^p\big(\varphi(a_1,...,a_p,\ )\cdot a_{p+1}\big)\bullet_{\alpha}\alpha^p(x),
                                \end{eqnarray*}
where
  $$\big(\varphi(a_1,...,a_p,\ )\cdot a_{p+1}\big)\bullet_{\alpha}\alpha^p(x)=\dl\sum_{i=1}^{n-1}[\alpha^p(x_{p+1}^1),...,\varphi(a_1,...,a_p,x_{p+1}^i ),...,\alpha^p(x_{p+1}^{n-1}),\alpha^p(x)].$$\\for all $a_i\in \mathcal{L}(\mathcal{N})$, $x\in \mathcal{N}$.
\end{thm}
\begin{proof}Let $\varphi\in\mathcal{C}^p_\mathcal{N}(\mathcal{N},\mathcal{N})$ and $(a_1\cdots a_{p+1})\in\mathcal{L} $ where $a_j=x_{1}^j\otimes\cdots\otimes x_{n-1}^j.$\\
 Then $\Delta\varphi\in \mathcal{C}^{p+1}_\mathcal{L}(\mathcal{L},\mathcal{L})$ and using \eqref{Leibnizcohomo} we can  to write $d=d_1+d_2+d_3$, where

\begin{align}
d_1\varphi(a_1,\cdots ,a_{p+1})& =\sum_{k=1}^p(-1)^{k-1} \big[ \alpha^{p-1}(a_k),\varphi(a_1,\cdots,\widehat{a_k},\cdots,a_{p+1})\big]\nonumber\\
d_2\varphi(a_1,\cdots ,a_{p+1})\ & =(-1)^{p+1} \big[ \varphi(a_1\otimes\cdots\otimes a_p),\alpha^{p-1}(a_{p+1}) \big]  \nonumber\\
d_3\varphi(a_1,\cdots ,a_{p+1})\ & =\sum_{1\leq k<j}^{p+1}{(-1)^k \varphi(\alpha(a_1)\otimes \cdots \otimes\widehat{a_k}\otimes \cdots\otimes\alpha(a_{j-1})\otimes [a_k,a_j]\otimes\alpha(a_{j+1})\otimes \cdots\otimes\alpha(a_{p+1}))}\nonumber
\end{align}
By \eqref{defdelta} we have
\begin{align}
& d_1\circ\Delta\varphi(a_1,\cdots ,a_{p+1})\nonumber\\
 \ & =\sum_{k=1}^p(-1)^{k-1} \big[ \alpha^{p-1}(a_k),\Delta\varphi(a_1,\cdots,\widehat{a_k},\cdots,a_{p+1})\big]\nonumber\\
 \ & =\sum_{k=1}^p(-1)^{k-1}\sum_{i=1}^{n-1} \big[ \alpha^{p-1}(a_k),\alpha^{p-1}(x_{p+1}^1)\otimes\cdots\otimes\varphi(a_1,\cdots,\widehat{a_k},\cdots,x_{p+1}^i)
 \otimes\cdots\otimes\alpha^{p-1}(x_{p+1}^{n-1})\big]\nonumber\\
 \ & =\sum_{k=1}^p(-1)^{k-1}\sum_{i>j}^{n-1}  \alpha^{p}(x_{p+1}^1)\otimes\cdots\otimes L(\alpha^{p-1}(x_k)).\alpha^{p-1}(x_{p+1}^j)\otimes\cdots\otimes\varphi(a_1,\cdots,\widehat{a_k},\cdots,x_{p+1}^i)
 \otimes\cdots\otimes\alpha^{p}(x_{p+1}^{n-1})\nonumber\\
 \ & +\sum_{k=1}^p(-1)^{k-1}\sum_{j>i}^{n-1}  \alpha^{p}(x_{p+1}^1)\otimes\cdots\otimes \varphi(a_1,\cdots,\widehat{a_k},\cdots,x_{p+1}^i)
 \otimes\cdots \otimes L(\alpha^{p-1}(x_k)).\alpha^{p-1}(x_{p+1}^j)\otimes\cdots\otimes\alpha^{p}(x_{p+1}^{n-1})\nonumber\\\ & +\sum_{k=1}^p(-1)^{k-1}\sum_{i=1}^{n-1}  \alpha^{p}(x_{p+1}^1)\otimes\cdots\otimes L(\alpha^{p-1}(x_k)).\varphi(a_1,\cdots,\widehat{a_k},\cdots,x_{p+1}^i)
 \otimes\cdots\otimes\alpha^{p}(x_{p+1}^{n-1})\nonumber\\
  \ & =\sum_{k=1}^p(-1)^{k-1}\sum_{i>j}^{n-1}  \alpha^{p}(x_{p+1}^1)\otimes\cdots\otimes L(\alpha^{p-1}(x_k)).\alpha^{p-1}(x_{p+1}^j)\otimes\cdots\otimes\varphi(a_1,\cdots,\widehat{a_k},\cdots,x_{p+1}^i)
 \otimes\cdots\otimes\alpha^{p}(x_{p+1}^{n-1})\nonumber\\
 \ & +\sum_{k=1}^p(-1)^{k-1}\sum_{j>i}^{n-1}  \alpha^{p}(x_{p+1}^1)\otimes\cdots\otimes \varphi(a_1,\cdots,\widehat{a_k},\cdots,x_{p+1}^i)
 \otimes\cdots \otimes L(\alpha^{p-1}(x_k)).\alpha^{p-1}(x_{p+1}^j)\otimes\cdots\otimes\alpha^{p}(x_{p+1}^{n-1})\nonumber\\
 \ & +\Delta\circ\delta_3\circ\varphi(a_1,\cdots ,a_{p+1})\nonumber\\
 \ & =\Lambda_1+\Lambda_2+\Delta\circ\delta_3\circ\varphi(a_1,\cdots ,a_{p+1})\nonumber
\end{align}
where
\begin{align}
 & \Lambda_1=\sum_{k=1}^p(-1)^{k-1}\sum_{i>j}^{n-1}  \alpha^{p}(x_{p+1}^1)\otimes\cdots\otimes L(\alpha^{p-1}(x_k)).\alpha^{p-1}(x_{p+1}^j)\otimes\cdots\otimes\varphi(a_1,\cdots,\widehat{a_k},\cdots,x_{p+1}^i)
 \otimes\cdots\otimes\alpha^{p}(x_{p+1}^{n-1})\nonumber\\
 \ & \Lambda_2=\sum_{k=1}^p(-1)^{k-1}\sum_{j>i}^{n-1}  \alpha^{p}(x_{p+1}^1)\otimes\cdots\otimes \varphi(a_1,\cdots,\widehat{a_k},\cdots,x_{p+1}^i)
 \otimes\cdots \otimes L(\alpha^{p-1}(x_k)).\alpha^{p-1}(x_{p+1}^j)\otimes\cdots\otimes\alpha^{p}(x_{p+1}^{n-1})\nonumber
\end{align}
Similarly we can prove that
\begin{align}
 d_2\circ\Delta\varphi(a_1,\cdots ,a_{p+1})
 \ & =\Delta\circ\delta_4\varphi(a_1,\cdots ,a_{p+1})\nonumber
\end{align}
and
\begin{align}
& d_3\Delta\circ\varphi(a_1,\cdots ,a_{p+1})\nonumber\\
 \ & =\sum_{1\leq k<j}^{p}{(-1)^k \Delta\circ\varphi(\alpha(a_1)\otimes \cdots \otimes\widehat{a_k}\otimes \cdots\otimes\alpha(a_{j-1})\otimes [a_k,a_j]\otimes\alpha(a_{j+1})\otimes \cdots\otimes\alpha(a_{p+1}))}\nonumber\\
\ & +\sum_{ k=1}^{p+1}{(-1)^k \varphi(\alpha(a_1)\otimes \cdots \otimes\widehat{a_k}\otimes \cdots\otimes\alpha(a_{p})\otimes [a_k,a_{p+1}])}\nonumber\\
\ & =\Delta\circ\delta_1\varphi(a_1,\cdots ,a_{p+1})+\Delta\circ\delta_2\varphi(a_1,\cdots ,a_{p+1})\nonumber\\
 \ & +\Lambda'_1+\Lambda'_2\nonumber
\end{align}
where $\Lambda'_1=-\Lambda_1$ and $\Lambda'_2=-\Lambda_2$.\\

Finally we have
$$d\circ\Delta=d_1\circ\Delta+d_2\circ\Delta+d_3\circ\Delta=\Delta\circ\delta_3+\Delta\circ\delta_4+
\Delta\circ\delta_1+\Delta\circ\delta_2=\Delta\circ\delta$$
where $\delta=\delta_1+\delta_2+\delta_3+\delta_4$ as defined in Proof  \ref{opercob}.
\end{proof}
\begin{remark}If $d^2=0$, then $\delta^2=0$.\\
In fact,we have $d\circ\Delta=\Delta\circ\delta$, then $$\Delta\circ\delta^2 =\Delta\circ\delta\circ\delta =d\circ\Delta\circ\delta=d\circ. d\circ\Delta=d^2\circ\Delta=0.$$
\end{remark}

\bibliographystyle{amsplain}
\providecommand{\bysame}{\leavevmode\hbox to3em{\hrulefill}\thinspace}
\providecommand{\MR}{\relax\ifhmode\unskip\space\fi MR }
\providecommand{\MRhref}[2]{%
  \href{http://www.ams.org/mathscinet-getitem?mr=#1}{#2}
}
\providecommand{\href}[2]{#2}

\end{document}